\documentclass[letterpaper, 10 pt, conference]{ieeeconf}  
\usepackage[T1]{fontenc}

\IEEEoverridecommandlockouts                              
\usepackage{natbib} 
\usepackage{amsmath, amssymb}
\usepackage{algpseudocode,algorithm,algorithmicx} 
\usepackage{cancel} 
\usepackage{graphicx}
\usepackage{tabularx}
\usepackage{booktabs}
\usepackage[symbol]{footmisc}
\usepackage{siunitx}
\usepackage{hyperref}
\usepackage{flushend} 
\usepackage[dvipsnames]{xcolor}
\usepackage{cleveref} 
\usepackage{makecell} 
\usepackage{soul} 
\usepackage{wrapfig} 

\title{\LARGE \bf
Realistic Extreme Behavior Generation for Improved AV Testing
}

\author{%
Robert Dyro\textsuperscript{1,4}, %
Matthew Foutter\textsuperscript{2,*}, %
Ruolin Li\textsuperscript{1}, %
Luigi Di Lillo\textsuperscript{3,5},
Edward Schmerling\textsuperscript{4}, \\
Xilin Zhou\textsuperscript{3},  %
Marco Pavone\textsuperscript{1,4}
\thanks{%
\textsuperscript{*} Corresponding Author. %
\textsuperscript{1} Dept. of Aeronautics and Astronautics, Stanford University.
\textsuperscript{2} Dept. of Mechanical Engineering, Stanford University.
\textsuperscript{3} Swiss Reinsurance Company Ltd. %
\textsuperscript{4} NVIDIA Corp. %
\textsuperscript{5} Autonomous Systems Lab, Stanford University (Research Affiliate). Contact: %
\texttt{\{rdyro, mfoutter, ruolinli, pavone\}@stanford.edu}, \texttt{eschmerling@nvidia.com}, \texttt{\{xilin\_zhou, luigi\_dilillo\}@swissre.com} %
}
}

\DeclareMathOperator{\qr}{\tt{qr}}
\DeclareMathOperator{\eig}{\tt{eig}}
\DeclareMathOperator{\minimize}{minimize}
\DeclareMathOperator{\suchthat}{such~that}
\DeclareMathOperator{\argmin}{argmin}

\newcommand\norm[1]{\ensuremath {\left|\left|{#1}\right|\right|}}

\newcommand\myparagraph[1]{\noindent\textbf{#1}}

\renewcommand{\baselinestretch}{0.938}
\begin{document}

\maketitle
\thispagestyle{empty}
\pagestyle{empty}

\begin{abstract}
This work introduces a framework to diagnose the strengths and shortcomings of Autonomous Vehicle (AV) collision avoidance technology with synthetic yet \emph{realistic} potential collision scenarios adapted from real-world, collision-free data. Our framework generates counterfactual collisions with diverse crash properties, e.g., crash angle and velocity, between an adversary and a target vehicle by adding perturbations to the adversary's predicted trajectory from a learned AV behavior model. Our main contribution is to ground these adversarial perturbations in realistic behavior as defined through the lens of data-alignment in the behavior model's parameter space. Then, we cluster these synthetic counterfactuals to identify plausible and representative collision scenarios to form the basis of a test suite for downstream AV system evaluation. We demonstrate our framework using two state-of-the-art behavior prediction models as sources of realistic adversarial perturbations, and show that our scenario clustering evokes interpretable failure modes from a baseline AV policy under evaluation.


\end{abstract}

\section{Introduction}

Autonomous Vehicles (AVs) promise to increase the efficiency and safety of transportation without the need for a human operator. 
However, challenges in assessing and validating the performance of AVs in the presence of other road users, and the ``long tail'' of behaviors these other agents may exhibit around the AV, make this promise of improved safety presently elusive to realize.
One might consider using failure modes observed during deployment to iteratively improve the autonomy stack in the flavor of continual learning \cite{survey-ContinualLearning}, however, 1) the large majority of mature AV deployment data is often mundane and without failure; 2) the frequency of observed failures diminishes as the vehicle's safety is further improved; and 3) we wish to \emph{preemptively} avoid unsafe behavior that may cause serious harm to the occupant(s) of the vehicle and surrounding environment.

Therefore, a fundamental challenge in AV safety is forecasting unseen, difficult scenarios that may rarely arise in nominal deployment. While real-world data may be insufficient, AV safety testing is fortunately well positioned to use simulation to synthesize and test these potentially problematic scenarios.

\begin{figure}[t]
    \centering
    \includegraphics[width=0.9\linewidth]{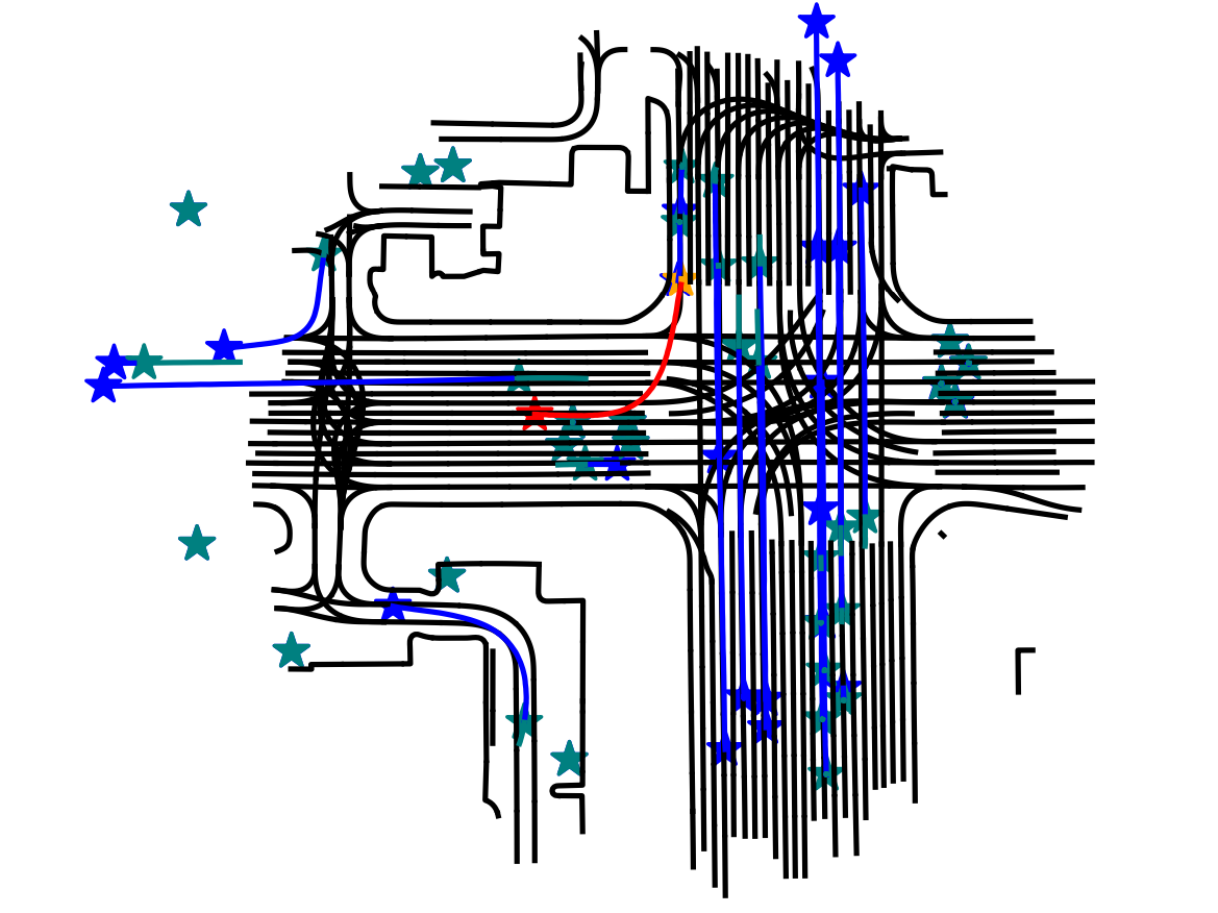}
    \caption{An example counterfactual collision from our framework generated by modifying a reference scenario from the Waymo Open Dataset \cite{waymo-open-dataset} with a perturbation derived from the MTR behavior model \cite{shi2022motion}. Our framework modifies the trajectory of an adversarial agent, using \emph{realistic} behavior perturbations, to encourage a collision with a target agent along a reference trajectory. In this example, all reference trajectories are highlighted in blue, while the adversary's trajectory is colored red. The non-adversarial agents' initial and final positions are highlighted by a green and blue star, respectively; the adversary's initial and final position are highlighted by a yellow and red star, respectively. We present a counterfactual in which the adversary performs an aggressive, over-wide right turn and collides with traffic stopped in the oncoming lane. Our approach offers valuable counterfactual scenarios -- grounded in the notion of realism -- on which to evaluate the maturity of existing AV collision avoidance technology.}
    \label{fig:mtr-collision}
    \vspace*{-0.6cm}
\end{figure} 

For example, using an efficient and sufficient representation like Bird's Eye View \cite{liu2021bevnet}, computer-generated scenarios and simulation-tested AV systems readily transfer learned safety behaviors to the real world in which case difficult scenario synthesis has the potential to increase AV robustness. This scenario generation and refinement must be automatic such that the generated features are more diverse and exhaustive as compared to manual generation.

Existing methods for generating adversarial scenarios face challenges in controlling realism and quality. While naive automated systems generate numerous scenarios, many are uninteresting. Therefore, selecting the most informative scenarios is crucial.
We combine adversarial optimization with a learned behavior model quantifying scenario likelihood to respectively achieve critical and realistic scenario generation.
Two key contributions are proposed:

\begin{enumerate}
\item \textbf{Realistic Scenario Generation Framework:} We present a general framework that optimizes perturbations to the weights of a pre-trained AV behavior model (which need only return a maximum a-posteriori output) to generate realistic counterfactual collision scenarios. We explicitly quantify the likelihood of the model's parameterization, with respect to its training data, to enforce that any perturbations to a vehicle's behavior are grounded in the realism captured by everyday interactions.

\item \textbf{Representative Scenario Clustering for Testing:} We cluster the synthetic counterfactual scenarios by diverse crash properties, e.g., crash angle and velocity, to reveal representative collision scenarios for downstream AV system testing. This counterfactual database offers a valuable test suite for developers, insurance companies and policymakers to evaluate the maturity and the risk of AV technology under diverse and realistic stress conditions.
\end{enumerate}

We demonstrate our contributions on two state-of-the-art behavior prediction models and evaluate our counterfactual database on a baseline AV policy to produce interpretable failure modes.

\section{Related Work}

\myparagraph{AV Behavioral Models:}\label{par:related-work-driving}
Behavior modeling is crucial for AV systems, enabling the prediction of other traffic participants' movements. Early works used simple dynamics models \cite{kong2015kinematic} and rule-based simulators \cite{dosovitskiy2017carla}. Neural-network-based models have since demonstrated superior performance in predicting future movement \cite{bansal2018chauffeurnet,cui2018multimodal}. These models often employ a Bird's Eye-View scene representation \cite{cui2018multimodal,varadarajan2021multipath} and encode the map \cite{shi2022motion} and the graph-interaction of participants \cite{ivanovic2018the,salzmann2020trajectron}. While some models output multiple possible futures \cite{cui2018multimodal,varadarajan2021multipath,shi2022motion}, others utilize a latent space approach for sampling \cite{ivanovic2018the,salzmann2020trajectron,rempe2021generating}. In this work, we desire a methodology to generate counterfactual scenarios that is compatible with all aforementioned architectures and is robust to future architecture development. Therefore, we only require the ability to evaluate a-posteriori likelihood from the behavior model, and are agnostic to specific input or internal representations.

\myparagraph{Counterfactual Scenario Generation:}
Generating unsafe, critical scenarios often involves adversarial approaches, either gradient-free, e.g., evolutionary algorithms \cite{biethahn1997evolutionary}, or gradient-based \cite{rempe2021generating,rempe2023trace,cao2022advdo}. Some works utilize photo-realistic simulators \cite{dosovitskiy2017carla} to generate plausible scenarios with 3D scene worlds \cite{okelly2018scalable}, while many recent works focus on simpler scene representations. This work adopts a simpler approach and gradient-based optimization, unlike \cite{wang2021advsim,vemprala2020adversarial,klischat2019generating}, with the added constraint that the adversarial perturbations to behavior are grounded in the \emph{realism} captured by a mature behavior model. Recent works have used generative models for scene construction \cite{xuNonediffscene,rempe2023trace}; however, these works do not prioritize robust realism for the generated trajectory. Unlike these works, we offer a definition of realism through the lense of data-alignment in the model's own parameter space.

\myparagraph{Clustering-based Scenario Exposition:}
Representative scenarios, characterized by static factors, e.g., weather and road geometry, and dynamic factors, e.g., trajectories and velocities, reveal varying levels of risk across different behavioral models. Such scenarios are typically derived from databases rich in crash data, employing various clustering methods. These databases may include real-world data from crash reports \cite{otte2003scientific,nitsche2017pre} or synthetic datasets designed for specific crash scenarios. Clustering techniques employed include KNN \cite{macqueen1967some,lloyd1982least}, k-medoids \cite{kaufman1990partitioning}, hierarchical clustering \cite{kaufman2009finding}, density-based clustering \cite{ester1996density}, and deep clustering methods \cite{guo2017deep}. In this work, we adopt the KNN technique to cluster synthetic counterfactual data by crash conditions, e.g., crash type and speed, for downstream AV safety evaluation. Therefore, our methodology offers a framework to reveal opportunities for further enhancements and areas in need of technological improvement within modern AV stacks.


\myparagraph{Probabilistic Learned-Model Analysis:}
To develop a general likelihood quantification method for deep behavior models, we turn to the probabilistic analysis of deterministic, learned maximum a-posterior models. Unlike more experimental methods \cite{kristiadi2020learnable, liu2021a,khan2019approximate,franchi2023make}, we utilize the general and model-agnostic Laplace approximation. While the Laplace approximation is generally intractable, efficient approximation methods exist \cite{ritter2018a,daxberger2021laplace}. We opt for a sketching-based approach \cite{tropp2016practical} which benefits from a strong theoretical analysis.

\section{Problem Formulation}

\subsection{Problem Setup}

In this work, we decompose driving scenarios into (i) a static scene description $S$ with semantic maps to identify the road and non-drivable area, e.g., road lanes, sidewalks, etc., (ii) the trajectories of $N$ non-adversarial vehicles $X = \{\textbf{X}^i\}_{i=1}^N$, sequences of 2D positions from a Bird's Eye View, and (iii) the trajectory of the adversarial vehicle $\textbf{X}^{\text{adv}}$. Each trajectory $\textbf{X}^i = [x_1^i, x_2^i, \dots, x_T^i]' \in\mathbb{R}^{T\times 2} \text{, where } T \in \mathbb{Z}_{+}$ is the final time, enumerates the state, $x_t^i$, for the $i$th agent at each time $t$; hence, $x^{\text{adv}}_1$ holds the initial condition for the adversarial agent. We assume access to a ground-truth trajectory for all non-adversarial agents in a real, reference scenario in which no collision occurs. The reference trajectory for all $N$ agents is given by $X_{\text{ref.}} = \{\textbf{X}^i_{\text{ref}}\}^N_{i=1}$.
Further, we assume access to a learned behavior model $f_{\text{bhv.}}(\theta, x^{\text{adv}}_1, X_{\text{ref.}}, S) =\textbf{X}_{\text{bhv.}}^{\text{adv}}$, parameterized by a vector of model weights $\theta\in\mathbb{R}^n$, which jointly processes the adversary's initial state, the reference trajectory of the $N$ non-adversarial agents and the scene description to produce a predicted trajectory for the adversary. Given these considerations, we desire a principled framework to choose the model weights $\theta$ such that the adversary's predicted trajectory from $f_{\text{bhv.}}$ 1) creates a critical, unsafe scenario through a collision between itself and another agent and 2) is realistic. In particular, this work's main contribution is the choice of model weights $\theta$ such that the adversary's resultant trajectory is grounded in the notion of \emph{realism}. Further, we seek a methodology to cluster these counterfactual scenarios according to similar descriptive features, e.g., velocity, thereby revealing the characteristics of the synthetic crashes to be used for safety evaluation.

\subsection{Objectives}

\begin{wrapfigure}{r}{0.2\textwidth}
    \centering
    \includegraphics[width=\linewidth]{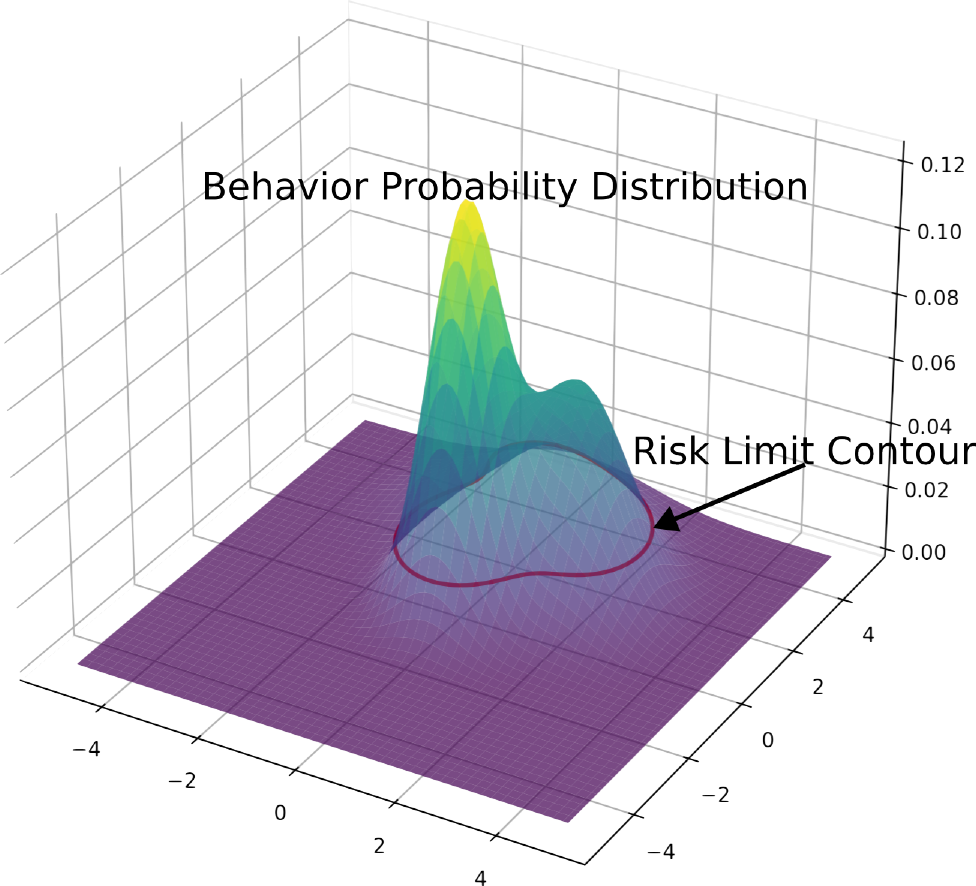}
    \caption{The risk contour $c$ offers a segmentation between realistic and unrealistic behavior with respect to a generic feed-forward behavior model.}
    \label{fig:risk-contour}
\end{wrapfigure}

The first objective is to synthetically generate realistic, counterfactual scenarios. Existing work on scenario generation, like \cite{rempe2021generating}, generates counterfactual scenarios by perturbing the behavior of existing scenarios. This process involves replacing the reference trajectory of an agent with the prediction from a behavior model such that we can adjust the model's weights to cause a synthetic collision. We aim to extend this generation framework by viewing realism through the lense of an envelope of larger than a constant $c$ behavior probability density as defined by the behavior model itself. We visually depict such a condition in \Cref{fig:risk-contour}.

As noted in \Cref{par:related-work-driving}, the challenge in incorporating a variety of behavior model architectures into a generation framework is that each model may have a different output format: To remain compatible with as many models as possible, and to guard against future model development, we seek a methodology that only utilizes a model's maximum a posteriori prediction to facilitate scenario generation.

With a counterfactual database, the second objective is identifying representative crash scenarios. We define representative scenarios by grouping crashes based on similarity and analyzing patterns within the resultant clusters, which involves three stages: (1) in \emph{crash reconstruction}, the crashes are reconstructed using core features, such as velocities and angles, to represent the severity of the collision; (2) in \emph{clustering}, the reconstructed crashes are clustered using an appropriate algorithm; and (3) in \emph{cluster analysis}, the generated clusters are analyzed using descriptive features to expose the characteristics influencing vehicle response.

This process aims to identify representative scenarios, i.e., group synthetic counterfactuals by similar properties, which can then be used to evaluate AV stacks and provide insights for improvement. Feature selection and clustering algorithms should reflect meaningful clusters, i.e., a lower bound on the minimum cluster size, and consistent cluster representations, i.e., a high Silhouette score, in order to guide future development.

\section{Approach}



Our proposed technique acknowledges the challenge of objectively defining and measuring realism, particularly in complex driving scenarios. While achieving perfect realism might be elusive, existing behavior prediction models, trained on massive datasets, implicitly capture a high degree of realism through their ability to forecast real-world driving behavior accurately. This observation suggests a potential path: if the optimal model weights represent the peak of realism within the training data, then deviations from this point can be used to define a quantifiable region of realism. We focus on perturbing the weights of a behavior prediction model within a neigborhood of magnitude $r\in\mathbb{R}_+$, a hyperparameter in our algorithm, from the optimal point. By analyzing the scaled distance of the deviation from optimality in the parameter space, we can establish a measurable metric that captures scenarios deemed realistic within the context of the training data. This technique offers a valuable solution, swapping the subjective notion of realism for a measurable notion of ``data-alignment," measured within the model's parameter space. 

In this section, we approximate the likelihood of a neural network's parameterization, in a scalable manner, using the Laplace approximation and techniques from matrix sketching. We also provide considerations to guide weight selection such that the likelihood of the model is preserved and define the adversarial objective we use to incentivize a collision. Then, we state our main contribution by formulating the space of realistic model parameters as a constraint set for adversarial optimization. Finally, we detail the descriptive features by which we cluster the counterfactual database and present a second, parameter-efficient matrix sketching technique leveraging Low-Rank Adaptation.

\subsection{Laplace Approximation}

We employ the Laplace approximation to principally approximate the likelihood levels of any maximum a-posteriori model. We begin with a second-order expansion of the loss function, $\ell$, parameterized by a vector of weights $\theta\in\mathbb{R}^n$, e.g., a neural network, about optimality, $\theta^{\star}$,
\begin{equation}
   \ell(\theta) \approx \ell(\theta^\star) + \cancelto{0}{\nabla\ell(\theta^\star)}^T \Delta \theta + \frac{1}{2} {\Delta \theta}^T H(\theta^\star) \Delta \theta,
\end{equation}
where $\Delta \theta = \theta - \theta^\star\in\mathbb{R}^n$ is a parameter perturbation from optimality.

Further, suppose that we define the loss function to correspond to the negative log posterior on the underlying weight distribution such that
\begin{equation}
\ell(\theta) = -\log p\left(\theta \mid \mathcal{D}\right) \approx -\log p\left( \theta^\star \mid \mathcal{D} \right) + \frac{1}{2} {\Delta \theta}^T H \Delta \theta,
\end{equation}
which leads to the interpretation of a Gaussian belief distribution over $\theta$, i.e.,
\begin{equation}
\theta \sim \mathcal{N}\left(\theta^\star, H^{-1}\right) ~~~~~~ \text{with} ~~ \Sigma = H^{-1}.
\end{equation}

However, the Hessian inverse covariance is computationally intractable as it exists in $\mathbb{R}^{n \times n}$ and $n$, the number of model parameters, is typically very large for deep learning models. Therefore, for likelihood estimation, we require an accurate and computationally efficient technique in order to develop an approximation of the Hessian inverse covariance.

\subsection{Symmetric Matrix Sketching}

We turn to numerical sketching to tractably approximate the Hessian's most important covariance components, which allows us to compress the top energy components of the Hessian efficiently. Matrix sketching is a technique that allows one to approximate a dense matrix $H$ as a product of low-rank factors, typically, $A Q$ where $A \in \mathbb{R}^{n \times k}$, $Q \in \mathbb{R}^{k \times n}$ and $k \ll n$. Alternatively, if $H$ is symmetric, we can use the 3-factor low-rank approximation $U D U^T$ where $U \in \mathbb{R}^{n \times k}$ and $D \in \mathbb{R}^{k \times k}$. Randomized sketching, an efficient sketching algorithm, is known to converge to a low-rank approximation of the top-energy components of the matrix $H$ \cite{tropp2016practical}. Therefore,  we can use such an algorithm to \emph{quantify the directions in the parameter space that lead to the fastest decrease in the likelihood.} We reproduce an efficient sketching algorithm from \cite{tropp2016practical} in the Appendix \Cref{sec:appendix-sketching}.


\subsection{Minimizing Probability Loss when Choosing \texorpdfstring{$\Delta z$}{D z}}

After establishing an estimate of the Hessian inverse covariance, we desire guidelines to reveal the weight perturbations, $\Delta\theta\in\mathbb{R}^n$, that result in a minimal decrease to the likelihood of the model's weights. Our objective is to develop a fixed-rank projection of an arbitrary vector $z \in \mathbb{R}^{|z|}$ defined as
\begin{equation}
\tilde{\Delta z} = (I - P P^T) \Delta z,
\end{equation}
such that $\tilde{\Delta z}^T H \tilde{\Delta z}$ is minimized. We formulate this objective as 

\begin{equation}
\begin{aligned}
\minimize_P & ~~~~ {\Delta z}^T (I - P P^T)^T H (I - P P^T) \Delta z  \\
\suchthat & ~~~~ P^T P = I \text{ and } \forall \Delta z \in \mathbb{R}^{|z|}.
\end{aligned}
\end{equation}
Then, this objective can be reformulated as 
\begin{equation}
\begin{aligned}
\minimize_P & ~~~~ \left|\left| H (I - P P^T) \right|\right|_\text{op} \\
\suchthat & ~~~~ P^T P = I,
\end{aligned}
\end{equation}
whose solution, for a fixed rank $k$, is to remove components from $z$ which correspond to the highest singular values in $H$, i.e.,
\[ P^\star = U[:, 1:k] ~~~ \text{where} ~~~ U S U^T = H. \]

\myparagraph{Projection Operation:} 
In linear algebra, a projection of the vector z $\in\mathbb{R}^n$ onto the range of a square matrix $A\in\mathbb{R}^{n\times n}$ produces the vector $\hat{z}\in\mathbb{R}^n$  given by
\[ \hat{z} = A (A^T A)^{-1} A^T z. \]
For an orthonormal matrix $P$, satisfying $P^T P = I$, this projection operation gives $\hat{z}$ as
\[ \hat{z} = P (\cancelto{I}{P^T P})^{-1} P^T z. \]
Therefore, we can write the complementary rejection of directions as
\[ \tilde{z} = z - \hat{z} = z - P P^T z = (I - P P^T) z, \]
in which case the operation $(I - PP^T)$ removes the components of $z$ in the range of $P$. 

\subsection{Adversarial, Collision-inducing Objective}

In order to incentivize a collision between two vehicles in the environment, we formulate an adversarial objective for the adversary in the scenario. From the problem formulation, we assume we have access to a behavior model $f_{\text{bhv.}}$, parameterized by weights $\theta \in\mathbb{R}^n$, in order to predict the adversary's state trajectory. We construct a simple collision loss with these primitives, akin to the formulation in \cite{rempe2021generating}, by choosing $\theta$ to minimize the minimum separation over time of the distance between the adversarial vehicle, $\textbf{X}^{\text{adv}}$, and each targeted neighbor, $\textbf{X}^{\text{target}}_{\text{ref.}},$ in the optimization batch with the added constraint that we restrict the search space of $\theta$ to the set of realistic model parameters, $\mathcal{C}_{\text{realism}}$:

\begin{equation}
\begin{aligned}
\minimize_{\theta} ~~~ & f_0(\theta) := \operatorname{softmin}_T \norm{\textbf{X}^\text{adv} - \textbf{X}^\text{target}_{\text{ref.}}}_2 &\\
\suchthat ~~~ & \textbf{X}^\text{adv} := f_\text{bhv.}(\theta, \textbf{X}_1^{\text{adv}}, X_\text{ref.}, S) \in \mathbb{R}^{T \times 2}, & \\
& S \equiv \text{scene description}, &\\
& \theta \in \mathcal{C}_\text{realism}. &
\end{aligned}
\end{equation}

\subsection{Optimization with Projection Feasibility Constraint}
Finally, we offer a definition for the space of realistic model parameters as a constraint set for the aforementioned adversarial objective, $f_0$. We turn to the projected gradient descent method in the presence of constraints. For the minimization of an adversarial loss $f_0$ over the domain of a model's parameters, $\theta\in\mathbb{R}^n$, we ensure realism by 1) confining $\theta$ to a neighborhood about $\theta^{\star}$, the pretrained optimum, of magnitude the hyperparameter $r \in \mathbb{R}_+$ using the L2 norm and 2) enforcing that weight perturbations must be orthogonal to the range of $P$ spanned by the directions leading to fastest decrease in the likelihood of the model; hence, weight perturbations used in the minimization of $f_0$ preserve likely behavior from the perspective of the behavior model:
\begin{equation}
\begin{aligned}
\minimize_\theta ~~~ & f_0(\theta) \\
\suchthat ~~~ & \mathcal{C}_\text{realism} := \begin{cases} 
\norm{\Delta \theta}_2 = \norm{\theta - \theta^\star}_2 \leq r \\
P P^T \Delta \theta = P P^T \left(\theta - \theta^\star\right) = 0
\end{cases}.
\end{aligned}
\end{equation}
Given this formulation above, the constraints projections are provided as
\[ \Pi_P (\theta) = 
\argmin_{P P^T (\tilde{\theta} - \theta^\star) = 0} \norm{\tilde{\theta} - \theta}_2^2 = \theta - P P^T (\theta - \theta^\star), \]
and 
\begin{equation}
\begin{aligned}
\Pi_{\norm{\Delta \theta} \leq r} 
& = \begin{cases}
\theta ~~ \text{if} ~~ \norm{\theta - \theta^\star}_2^2 \leq r, \\
\theta^\star + \frac{r}{\norm{\theta - \theta^\star}_2^2} (\theta - \theta^\star) ~~ \text{otherwise}
\end{cases}.
\end{aligned}
\end{equation}

\subsection{Clustering-based Representative Critical Scenarios}

After generating a counterfactual database, we aim to identify representative collision scenarios for testing autonomy stacks. We group similar crashes using clustering to extract these cases.

We focus on \emph{core} features related to vehicle dynamics before the crash, reflecting our key concerns about crashes. These features include velocity, relative velocity, and angle of impact. Additionally, we include features signifying the tested vehicle's response, such as response time. We categorize the angle feature into crash type, e.g., contrasting, side-left, side-right and chasing. We define a "contrasting" crash to occer when the adversary is in the oncoming lane, while a "side-left", or "side-right", crash to occur with the adversary to the left, or right, of the target; we define a "chasing" crash to occur when the adversary begins trailing the target.

For \emph{descriptive} features, we use: Adversarial vehicle speed $v_{\text{a}}$; Longitudinal relative speed $\Delta v_x$; Lateral relative speed $\Delta v_y$; Angle of impact $\gamma$; Crash type $\beta$; and response time $t_r$. \Cref{tab:descriptive-features} in the Appendix \Cref{sec:appendix-features} gives detailed descriptions of these features.

We use K-Nearest Neighbors (KNN) for clustering and evaluate the quality of the chosen clusters using the average Silhouette score \cite{shahapure2020cluster}. We impose a constraint on the smallest cluster size to prevent the formation of non-representative clusters. The choice of this size is dataset-specific and is discussed in the results section.

\subsection{Low-Rank Adaptation in Efficient Sketching}

Also, we propose an alternative to sketching the full Hessian inverse covariance: we sketch the Hessian in a lower dimensional, parameter-efficient space with Low-Rank Adaptation. For every NN layer representable as a matrix $W$, Low-Rank Adaptation enables parameter-efficient perturbations to the layer's weight space by adding a low-rank additive factor, i.e.,
\[
\tilde{W} = W + A B,
\]
where $A \in \mathbb{R}^{m \times k}$ and $B \in \mathbb{R}^{k \times n}$. For
$k \ll \{m, n\}$ and, importantly, at initialization,
\[
A_{ij} \sim \mathcal{N}(0, \sigma) ~~~~ B_{kl} = 0.
\]
In this approach, we sketch merely the factors $A$ and $B$ at their initialization point. In the Appendix \Cref{sec:appendix-gip}, we investigate the convergence of the Global Information Projection Matrix, which is the ability to extract the most energetic directions in the parameter space from the Laplace approximation for each sketching approach.



\section{Results}
In this section, we encapsulate the descriptive features from our synthetic counterfactual database by forming diverse, representative crash clusters, and we evaluate the response of a baseline AV policy on this dataset to expose the policy's strengths and weaknesses.

\begin{figure}[t]
    \centering
    \includegraphics[width=0.9\linewidth]{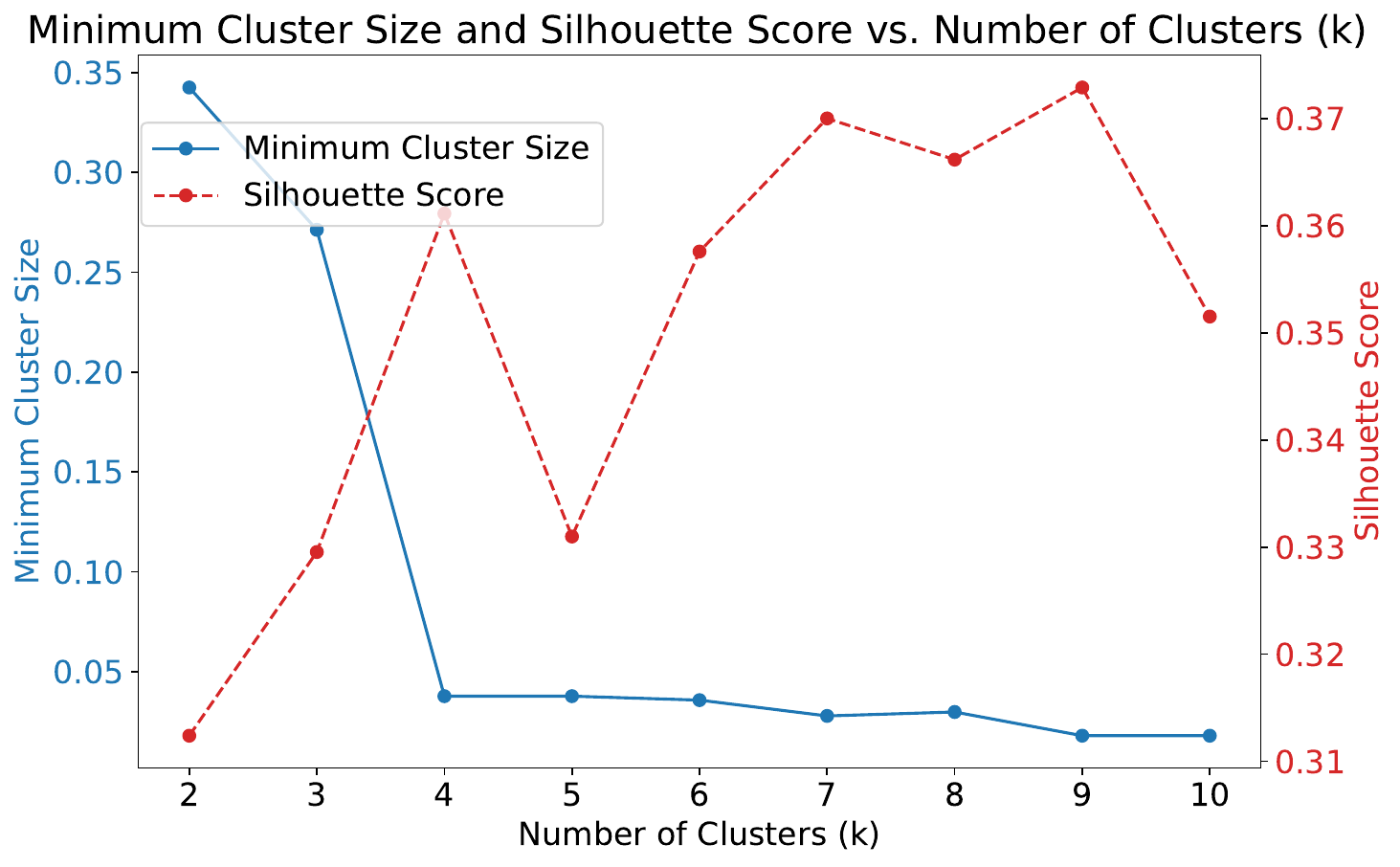}
    \caption{The number of clusters balances a representative set of scenarios, measured through the minimum cluster size, and the effectiveness of the clustering algorithm, measured by the Silhouette score. $k$, the number of clusters, achieves a low minimum cluster size and high Silhouette score for $k=8$.}
    \label{fig:k-choice-mtr}
\end{figure}

\begin{figure}[b]
    \centering
    \includegraphics[width=0.9\linewidth]{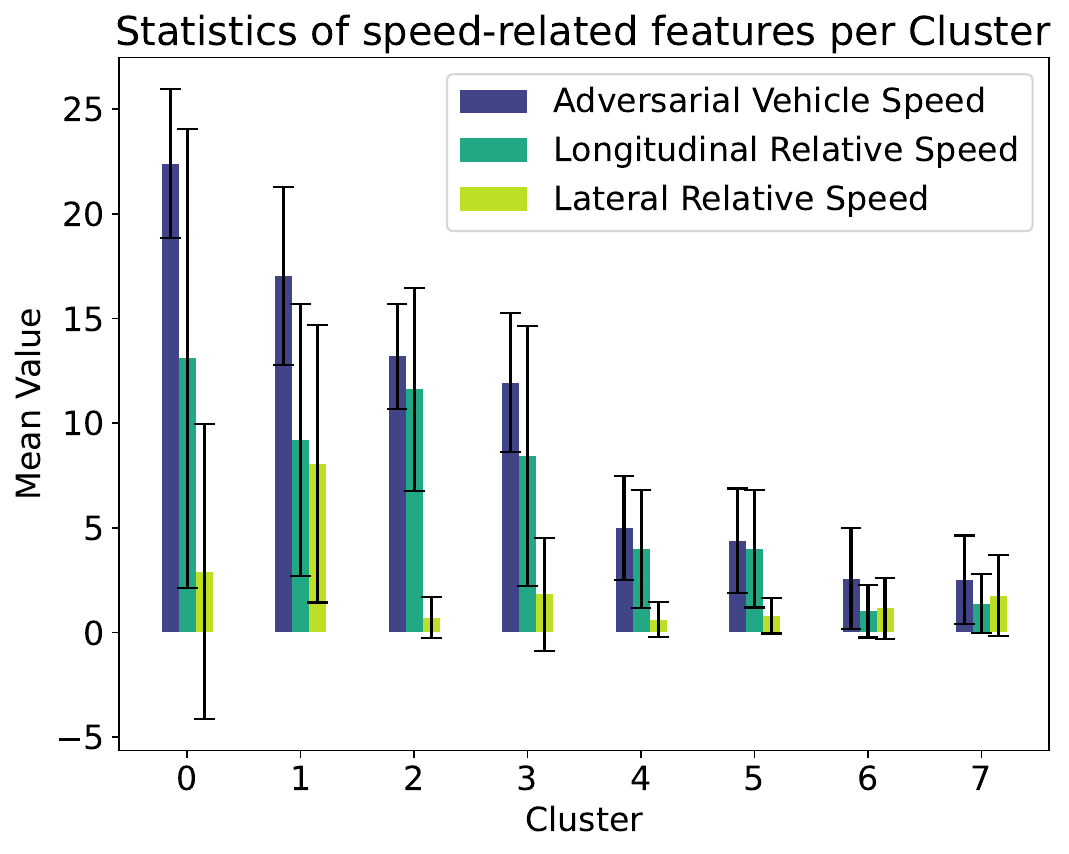}
    \caption{Clusters are numbered in descending order with respect to the severity of the collision by virtue of the mean adversarial vehicle speed.}
    \label{fig:cluster_des_srr_mtr}
\end{figure}
Our proposed counterfactual generation pipeline is agnostic to the underlying behavior model. We, therefore, use two distinct behavioral models: MTR and Trajectron++ \cite{shi2022motion, salzmann2020trajectron}. As discussed in \Cref{sec:behavioral_model-discussion}, we only present our findings for one representative behavioral model, i.e., MTR, in the main body. We evaluate the response of a baseline AV policy on the
target, non-adversarial vehicle with limited reactivity in each
counterfactual. The policy monitors the position and velocity of all agents, predicts the time and minimum distance to collision, assuming constant velocity extrapolation, and chooses to emergency brake if either measure is below a prescribed value; otherwise, the policy follows the reference trajectory using a single-step Model Predictive Control, obeying acceleration limits. The exact implementation of the baseline policy is provided in the Appendix \Cref{sec:appendix-baseline-policy}.

For the MTR behavior model, the collisions are generated with a prescribed realism hyperparameter $r = 0.003$ whose value is determined as discussed in \Cref{sec:realism-discussion}. The generated counterfactual dataset comprises a total of 505 collisions. We have set a threshold for the smallest cluster size at 3\% of the total dataset, aiming to identify clusters of relevant and representative scenarios. Therefore, we chose $k=8$ clusters in this case as shown in \Cref{fig:k-choice-mtr}. Final detailed statistics of the clusters are presented in \Cref{tab:statistics_mtr}. The clusters are labeled in descending order of the adversarial vehicle speed, $v_a$, from more intense to less intense testing conditions, i.e., see \Cref{fig:cluster_des_srr_mtr}.

\begin{table*}[ht]
    \centering
    \scriptsize
    \begin{tabularx}{\textwidth}{l *{8}{>{\centering\arraybackslash}X}}
        \toprule
        Descriptive Feature & Cluster 0 & Cluster 1 & Cluster 2 & Cluster 3 & Cluster 4 & Cluster 5 & Cluster 6 & Cluster 7 \\
        \midrule
        \textbf{Adversarial Vehicle Speed $v_a$} & 22 (3.5) & 17 (4.2) & 13 (2.5) & 12 (3.3) & 5.0 (2.5) & 4.4 (2.5) & 2.6 (2.4) & 2.5 (2.1) \\ 
        \textbf{Longitudinal Relative Speed $\Delta v_x$} & 13 (11) & 9.2 (6.5) & 12 (4.8) & 8.4 (6.2) & 4.0 (2.8) & 4.0 (2.8) & 1.0 (1.3) & 1.4 (1.4) \\ 
        \textbf{Lateral Relative Speed $\Delta v_y$} & 2.9 (7.1) & 8.1 (6.6) & 0.70 (0.98) & 1.8 (2.7) & 0.61 (0.83) & 0.79 (0.84) & 1.2 (1.5) & 1.8 (1.9) \\ 
        \midrule
        \textbf{Crash Type $\beta$} \\
        Chasing & \textbf{12 [43\%]} & 4 [27\%] & \textbf{91 [100\%]} & 0 & \textbf{120 [98\%]} & 0 & 0  & 14 [17\%] \\ 
        Contrasting & 7 [25\%] & 3 [20\%] & 0 & \textbf{51 [91\%]} & 0 & \textbf{67 [100\%]} & 10 [23\%] & 0 \\ 
        Side-left & 4 [14\%] & \textbf{7 [47\%]} & 0 & 1 [2\%] & 3 [2\%] & 0 & \textbf{34 [77\%]} & 0 \\ 
        Side-right & 5 [18\%] & 1 [7\%] & 0 & 4 [7\%] & 0 & 0 & 0 & \textbf{69 [83\%]} \\ 
        \midrule
        \textbf{Tested Vehicle Response} \\
        No response & 0 & 0 & 0 & 0 & 9 [7\%] & 12 [18\%] & 26 [59\%] & 46 [55\%] \\
        Response & 28 [100\%] & 15 [100\%] & 91 [100\%] & 56 [100\%] & 112 [93\%] & 55 [82\%] & 18 [41\%] & 37 [45\%] \\ 
        \midrule
        \textbf{Response Time $t_r$} & 6.7 (1.4) & 7.2 (0.68) & 5.2 (1.5) & 6.1 (1.8) & 5.0 (1.9) & 5.6 (1.7) & 6.0 (2.4) & 2.8 (2.6) \\
\bottomrule
\end{tabularx}
    \caption{Synthetic counterfactuals are clustered according to descriptive crash properties. We find that the baseline policy responds best to high-speed collisions and struggles to detect "side-left" and "side-right" crash types at low speeds.\label{tab:statistics_mtr}}
    \vspace*{-0.5cm}
\end{table*}

These diverse clusters reveal distinct patterns in the behavior of tested vehicles under various crash scenarios, providing insights into the operational challenges and efficiencies of a candidate AV policy.

\textbf{High-Speed Clusters:} Clusters 0 and 1, characterized by high testing speeds, i.e., 22.41 m/s and 13.19 m/s, respectively, predominantly feature "chasing" and "side-left" crashes. We assume the tested vehicle to have produced a "response" if the target correctly reacts to an imminent collision with the adversary by beginning to brake. Therefore, the 100\% response rate in both clusters 0 and 1 suggests the effective handling of high-stress conditions, indicating advanced detection and robust response from the baseline policy.

\textbf{Mid-Speed Clusters:} Clusters 2 and 3 demonstrate a lower mean testing speed, i.e., 13.19 m/s and 11.92 m/s, respectively, with a dominant number of "chasing" and "contrasting" crashes. We define a "contrasting" crash to occur when the adversary is in the oncoming lane. The 100\% response rate in both clusters highlights the system's effectiveness in detecting and reacting to both tailing and oncoming vehicle collisions.

\textbf{Low-Speed, High Complexity Clusters:} Cluster 4 showcases an almost complete response rate of 93\% despite a lower testing speed of 5.01 m/s, indicating high vigilance from the AV policy in tracking vehicles from behind at slower speeds. Cluster 5 reveals challenges in less predictable crash scenarios with an 18\% no-response rate to "contrasting" crashes, suggesting potential shortcomings in head-on detection at low speeds.

\textbf{The Most Challenging Clusters:} Clusters 6 and 7, with the lowest testing speeds, are dominated by "side-left" and "side-right" crashes, respectively, both presenting high no-response rates at 59\% and 55\% as shown in \Cref{tab:statistics_mtr}. This poor response rate may indicate difficulties for the policy to detect and react to lateral threats at low speeds, highlighting a potential weaknesses in the baseline's limited reactivity.



We offer \Cref{tab:statistics_mtr} to present a holistic view of our counterfactual database; however, any one counterfactual is only valuable to the extent the scenario is realistic. The patterns across clusters highlight that these counterfactuals demonstrate \emph{diverse} crash characteristics, e.g., a variety in crash severity and crash type, which, more importantly, are also grounded in the notion of \emph{realism} by virtue of the realistic behavior captured by the pre-trained behavior model. Therefore, this counterfactual generation pipeline offers the opportunity for critical insights to strengthen AV technology with a broad spectrum of real-world scenarios. For example, severe crashes with high response rates present opportunities for continued improvement, while less severe instances, with relatively high no-response rates, offer critical areas in need of technological improvements, particularly in sensor accuracy and response algorithms for the baseline policy.



\section{Discussion}

\subsection{The choice of \texorpdfstring{$r$}{r}}\label{sec:realism-discussion}
The realism radius $r$ is a hyperparameter trading off fidelity to the most data-aligned, i.e, ``realistic", prediction and aggressiveness of the desired behavior.
One cannot prescribe a single value $r$ that would work in all cases, but we make two observations: 1) by expressing $r$ as a constraint rather than a soft objective, it is not affected by scaling the adversarial objective function; and 2) because $r$ is always a scalar, we can use conformal prediction \cite{shafer2008tutorial, angelopoulos2023conformal} to determine its value in practice with a small calibration set. In our experiments, we set $r$ using the bisection method until each scenario contains a non-severe collision by visually inspecting the resultant behaviors on a small ($n=10$) calibration set.

\subsection{Selection of behavioral models}\label{sec:behavioral_model-discussion}

We highlight that our method is versatile and applicable to any parametric behavioral model. We demonstrate our pipeline on two representative behavior models: a regressive model \cite{salzmann2020trajectron} and a transformer-based model \cite{shi2022motion}. We have detailed the results for the transformer-based model \cite{shi2022motion} in the paper's main body, while, for brevity, we present the findings for the regressive model \cite{salzmann2020trajectron} in the Appendix \Cref{sec:appendix-trajectron++} with discussion. Further, we visually depict a synthetic counterfactual collision from each behavior model in \Cref{fig:mtr-collision}, on the main body's first page, and in \Cref{fig:generated-collisions-trajectron-1} in the Appendix \Cref{sec:appendix-collisions}.

\section{Conclusions}
This work proposes a novel framework for generating realistic and challenging scenarios for AV testing using limited collision-free data. The approach combines adversarial optimization with a trained behavior model, enabling the quantification of scenario likelihood and ensuring the generation of realistic counterfactual scenarios. Our framework leverages the Laplace approximation and sketching techniques for computationally scalable likelihood estimation, overcoming the limitation of previous work. The method is agnostic to the underlying behavior model or adversarial loss and is scalable, making it a versatile tool for testing AV systems in unsafe situations with a collision. We show a scaling-friendly parametrization based on the Low-Rank Adaptation reformulation. We demonstrate the effectiveness of our work on two state-of-the-art behavior prediction models and two distinct driving datasets. Furthermore, we identify representative and diverse crash conditions among the synthetic data based on KNN clustering for downstream policy evaluation. Our framework offers a systematic approach to generate \emph{realistic} counterfactual collisions with various AV behavior models; therefore, our approach offers diverse, high-fidelity and simulated stress conditions to reveal the strengthes and weaknesses of AV technologies for AV stakeholders, e.g., car manufacturers, insurance companies and government regulators.


\section*{Acknowledgment}

\noindent We thank Swiss Re for their support in conducting this work.
\newpage
\nocite{*}
\bibliography{references.bib}
\flushcolsend

\clearpage

\section*{APPENDIX}
\subsection{LoRA Hessian Computation}

We observe that despite reducing the required memory and computation by taking only the top eigenvector components of the inverse covariance matrix, $P\in\mathbb{R}^{n \times k}$, the representation can still be too large, some large ML models fit in memory only once, implying $k = 1$. and only one eigenvector component can be used. 

Instead, we turn to another compressed representation of the NN parameter representation. Importantly, this is not a third approximation in our algorithm: we perform the Laplace approximation in a smaller parameter space. Many choices exist for a reduced parameter space, the most obvious being a subset of network parameters. However, inspired by recent advances in fine-tuning large language models, we reparametrize the network using the Low-Rank Adaptation parameter space \cite{hu2021lora}.

Since $P$ is constructed, we can use a network \textbf{reparametrization} before constructing the Laplace approximation, i.e.,
\[
\forall W_i ~~~~ \text{let} ~~~~ \widetilde{W}_i = W_i + B_i A_i,
\]
where $B_i = 0$ and $A_i \sim \mathcal{N}(0, \sigma I)$.

For a new parametrization, for the Laplace approximation, we need $\nabla_p \ell =
0$ and $\nabla_p^2 \ell \neq 0$. Therefore, because
\[
\nabla_{\{A, B\}} \ell(W_i + B A) = \left(\nabla_{\widetilde{W}_i} \ell(\widetilde{W}_i) \right) \{ B, A^T \},
\]
and at trained network optimality $\nabla_{\widetilde{W}_i} \ell(\widetilde{W}_i) = 0$, then $\nabla_p \ell = 0$ and

\[ \nabla_A^2 \ell = (I \otimes B) H B = 0, \]
\[ \nabla_B^2 \ell = (I \otimes A^T) H A^T \neq 0. \]
Hence, we retain only $B_i$ parameters for the new parametrization.

\subsection{Sketching}\label{sec:appendix-sketching}

Below, in \Cref{alg:low-rank-sketch,alg:low-rank-sym-sketch,alg:low-rank-psd-sketching} we reproduce an efficient random sketching algorithms from \cite{tropp2016practical}. While \cite{tropp2016practical} offers several possible random sketching algorithms, we focus on the symmetric positive decomposition consisting of 2 factors
\[
    H \approx U S U^T.
\]

We specifically retrieve the positive semi-definite factor $S$ given that under Laplace method Gaussian approximation, the Hessian, or the inverse covariance matrix is positive semi-definite.

\subsubsection{Hessian matrix spectral decay in sketched matrices in \cite{salzmann2020trajectron, shi2022motion}}

\begin{figure}[t]
    \centering
    \includegraphics[width=\linewidth]{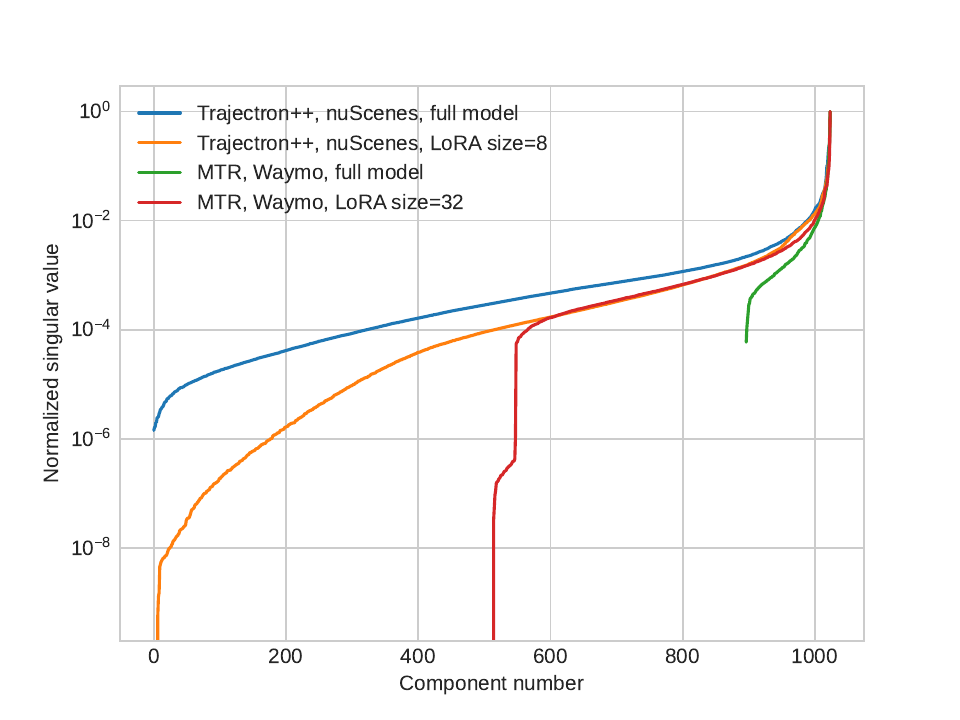}
    \caption{Singular values of the sketching decomposition of the Hessian
    matrix for the Trajectron++ and MTR models \cite{salzmann2020trajectron}\cite{shi2022motion}. Rapid decay of singular values indicates that a low-rank approximation captures most of the singular-value energy within the matrix.}
    \label{fig:singular-values}
\end{figure}



\subsubsection{Sketching algorithms applied in this work} Below, we concisely reproduce the sketching algorithms exploited in this work for sketching positive definite matrices comprising the covariance in the Laplace method as applied to the behavior model's parameter space. The sketching matrices are obtained as the result of \Cref{alg:low-rank-psd-sketching}, but \Cref{alg:low-rank-sketch,alg:low-rank-sym-sketch} are called as subroutines.

\vspace*{0.5cm}

\begin{algorithm}[H]
\begin{algorithmic}[1]
\Require{Left sketch matrix $\Psi$, Right sketch matrix $\Phi$}
\Require{Left sketch $W = \Psi A$, Right sketch $Y = A \Phi$}

\Function{\tt low\_rank}{$\Psi$, $\Phi$, $W$, $Y$}
\State $Q, \_ \leftarrow \qr(Y)$ \Comment{left orthonormal basis}
\State $U, T \leftarrow \qr(\Phi Q)$
\State $X \leftarrow T^{\dag} \left(U^* W\right)$ \Comment{triangular solve}
\State \Return $Q$, $X$ \Comment{$A \approx Q X$}
\EndFunction
\end{algorithmic}
\caption{Sketching Algorithms Reproduced From \cite{tropp2016practical}}
\label{alg:low-rank-sketch}
\end{algorithm}

\begin{algorithm}[H]
\begin{algorithmic}[1]
\Require{Left sketch matrix $\Psi$, Right sketch matrix $\Phi$}
\Require{Left sketch $W = \Psi A$, Right sketch $Y = A \Phi$}

\Function{\tt low\_rank\_sym}{$\Psi$, $\Phi$, $W$, $Y$}
\State $Q, X \leftarrow {\tt low\_rank}(\Psi, \Phi, W, Y)$
\State $U, T \leftarrow \qr([Q, X^*])$ \Comment{orthogonalize concatenation}
\State $T1 \leftarrow T[:, 1:k]$ and $T2 \leftarrow T[:, k+1:2 k]$ \Comment{split}
\State $S \leftarrow (T_1 T_2^* + T_2 T_1^*) / 2$ \Comment{symmetrize}
\State \Return $U, S$ \Comment{$A \approx U S U^*$}
\EndFunction
\end{algorithmic}
\caption{Sketching Algorithms Reproduced From \cite{tropp2016practical}}
\label{alg:low-rank-sym-sketch}
\end{algorithm}

\begin{algorithm}[H]
\begin{algorithmic}[1]
\Require{Left sketch matrix $\Psi$, Right sketch matrix $\Phi$}
\Require{Left sketch $W = \Psi A$, Right sketch $Y = A \Phi$}
\Function{\tt low\_rank\_psd}{$\Psi$, $\Phi$, $W$, $Y$}
\State $U, S \leftarrow {\tt low\_rank\_sym}(\Psi, \Phi, W, Y)$
\State $V, D \leftarrow \eig(S)$ \Comment{eigendecomposition}
\State $U \leftarrow U V$ \Comment{consolidate orthonormal factors}
\State $D \leftarrow \max(D, 0)$ \Comment{remove negative eigenvalues}
\State \Return U, D \Comment{$A \approx U D U^*$}
\EndFunction

\end{algorithmic}
\caption{Sketching Algorithms Reproduced From \cite{tropp2016practical}}
\label{alg:low-rank-psd-sketching}
\end{algorithm}

\subsection{Convergence of Global Information Estimate}\label{sec:appendix-gip}

\begin{figure}[b]
    \centering
    \includegraphics[width=\linewidth]{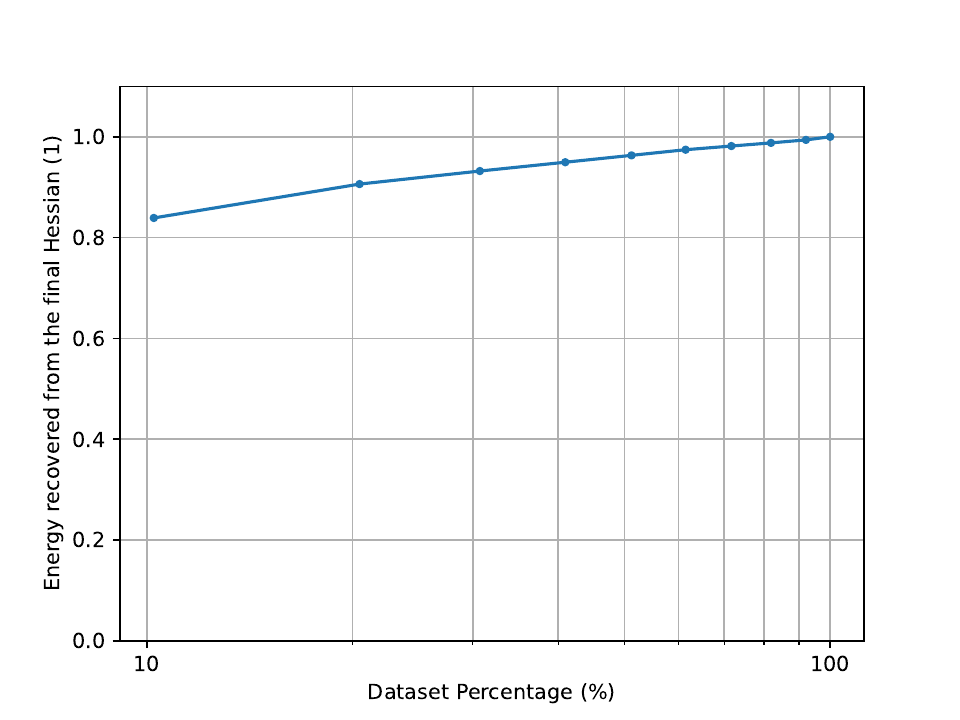}
    \caption{Global Information Projection (GIP) Matrix converges quickly for a
    fixed rank size $k$ with the number of scenarios in the training dataset.}
    \label{fig:hessian-convergence}
\end{figure}

We investigate the convergence of the Global Information Projection Matrix, which is the ability to extract the most energetic directions in the parameter space from the Laplace approximation. In \Cref{fig:hessian-convergence}, we show that the Hessian converges rapidly -- at 10\% of the dataset, the number of captured linear directions is above 80\%.

Likewise, for the Low-Rank Adaptation sketching case, we observe similar but more rapid convergence to the final extracted linear subspace in \Cref{fig:lora-hessian-convergence} presented in the Appendix. We theorize this is due to the inherently smaller parameter space through the use of Low-Rank Adaptation. The fast convergence in \Cref{fig:hessian-convergence} and \Cref{fig:lora-hessian-convergence} justifies that our fixed rank approximation of the Hessian only requires a fraction of the dataset; therefore, our proposed approach can be extended to large-scale behavior models trained on prohibitively large datasets. 

Additionally, we quantify the decay of the singular values in the symmetric positive-definite sketching decomposition of the Hessian matrix in \Cref{fig:singular-values} and observe, like other works, e.g., \cite{sharma2021sketching}, that the singular values of deep learned models decay rapidly. This result suggests that the sketching decomposition is an excellent approximation of the Hessian matrix.

\begin{figure}
    \centering
    \includegraphics[width=\linewidth]{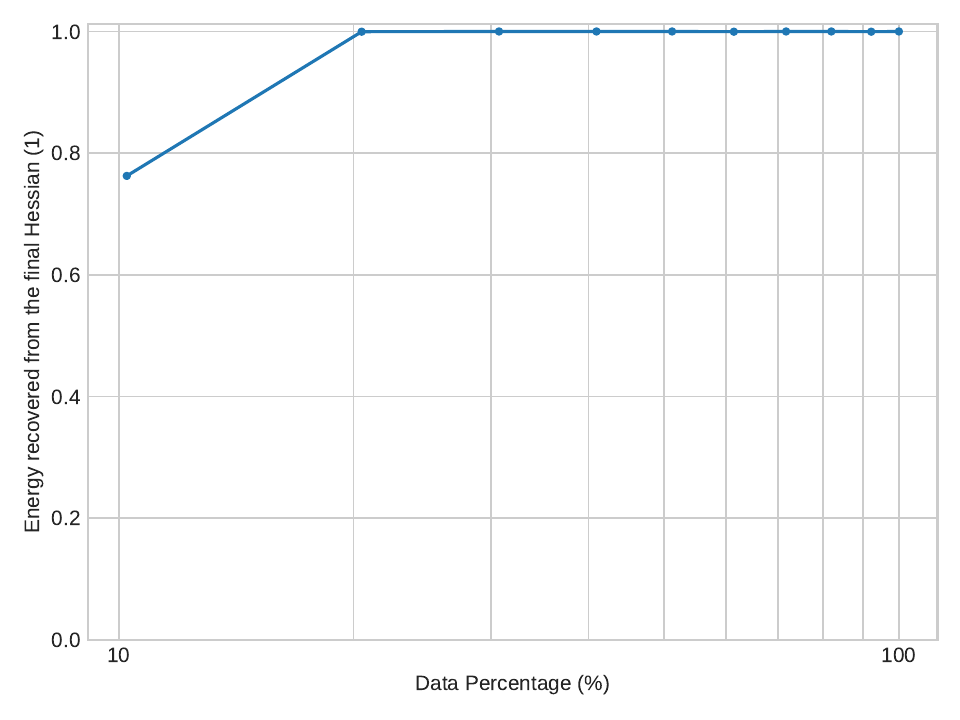}
    \caption{LoRA case: GIP matrix converges quickly for a fixed rank size $k$
    with the number of scenarios in the training dataset.}
    \label{fig:lora-hessian-convergence}
\end{figure}

\subsection{Baseline AV Policy For Evaluation}\label{sec:appendix-baseline-policy}

The baseline policy used in this work for modeling the reaction of other agents needed to be (i) adhering as closely to a reference trajectory (e.g., recorded behavior in the scene) but also (ii) reasonably realistic, allowing for an agent to react to the adversarial behavior directed at it. The policy monitors the position and velocity
of all agents, predicts the time and minimum distance to collision assuming constant velocity extrapolation, and chooses to 
emergency brake if either measure is below a prescribed value;
Otherwise, the policy follows the reference trajectory using a single-step Model Predictive Control, obeying acceleration limits. 

The mathematical formulation takes the form
\begin{align*}
\min_{a(t)} \quad & || x(t+1) - x_{ref}(t+1) ||^2 \\
& ~~~~~~~ + \frac{1}{10} || v(t+1) - v_{ref}(t+1) ||^2 \\
\text{s.t.} \quad & x(t+1) = x(t) + v(t) \Delta t, \\
& v(t+1) = v(t) + a(t) \Delta t, \\
& d(x(t+1), x_i(t+1)) \geq d_{min}, \\
& t_{coll.}(x(t+1), v(t+1), x_i(t+1), v_i(t+1)) \geq t_{min}, \\
& -a_{max} \leq a(t) \leq a_{max},
\end{align*}

\subsection{Collision Examples}\label{sec:appendix-collisions}

Examples of generated collisions are presented in the main body with \Cref{fig:mtr-collision} and in the Appendix with \Cref{fig:generated-collisions-trajectron-1}.

\begin{figure}[t]
    \centering
    \includegraphics[width=0.7\linewidth]{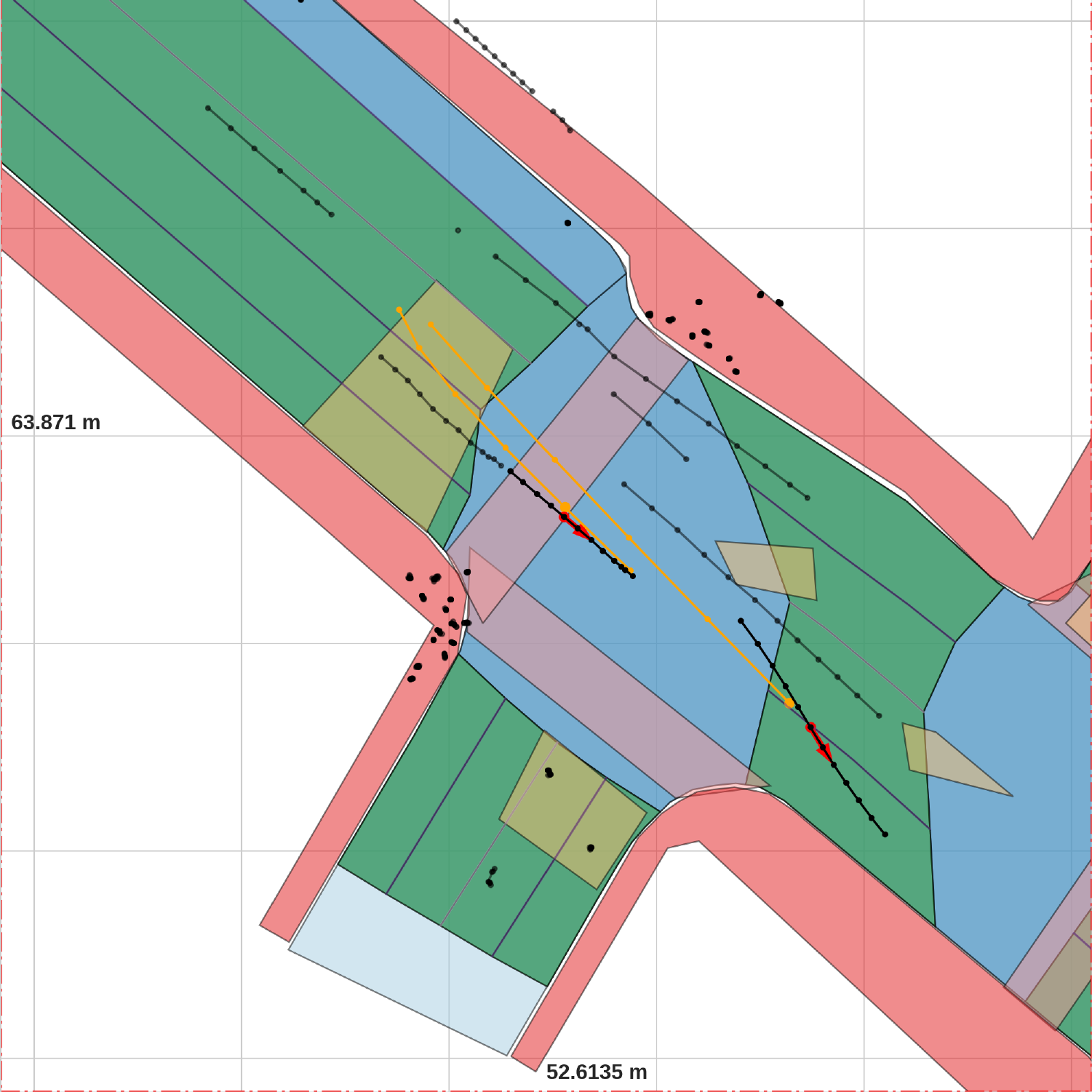}
    \caption{Example of generated collisions with other traffic participants using the nuScenes dataset and the Trajectron++ behavior model. The aggressive agent collides with the vehicle in the neighboring lane.}
    \label{fig:generated-collisions-trajectron-1}
\end{figure}


\begin{table*}
\small
\centering
\begin{tabular}{||p{4.5cm} | p{1.5cm} | p{8.5cm} ||} 
 \hline
Features &  Type & Description \\ [0.5ex] 
 \hline\hline
Adversarial vehicle speed $v_{\text{a}}$ &  Numeric & The magnitude of velocity for the adversarial vehicle right before the crash. \\
 \hline
Longitudinal relative speed $\Delta v_x$ &  Numeric & The maginitude of the relative velocity of the two crashing parties right before the crash in the longitudinal direction of the tested vehicle.\\
 \hline
Lateral relative speed $\Delta v_y$ &  Numeric & The magnitude of the relative velocity of the two crashing parties right before the crash in the lateral direction of the tested vehicle.\\
 \hline
Angle of impact $\gamma$ &  Numeric & The relative angle of two crashing parties at the crash. \\
 \hline
Crash type $\beta$ & Categorical &  Types of crashes categorized by the angle of impact. Types include contrasting, chasing, side-left and side-right.\\
 \hline
Response time $t_r$ & Numeric & The time, in seconds, before the tested vehicle responds, i.e., brakes, to the crash, if there is a response.\\
 \hline
\end{tabular}
\caption{\normalfont Core and descriptive features for generated collision clustering analysis.}
\label{tab:descriptive-features}
\end{table*}

\subsection{Quantifiable dynamics features used for behavioral clustering}\label{sec:appendix-features}

Our clustering analysis focuses on identifying representative critical scenarios for testing autonomous vehicle stacks by grouping similar crashes. To achieve this, we employ a combination of core and descriptive features. The core features, which reflect key concerns about crashes, include velocity, relative velocity, and angle of impact. Additionally, the response time of the tested vehicle is considered a core feature.

Descriptive features, on the other hand, categorize the angle feature into distinct crash types (contrasting, side-left, side-right, chasing). These features, along with the core features, are then used for clustering. The study utilizes the K-Nearest Neighbors (KNN) algorithm for clustering, ensuring representativeness by imposing a constraint on the smallest cluster size to prevent the formation of non-representative clusters. 

The detailed description of the features are contained in \Cref{tab:descriptive-features}.

\subsection{Representative Critical Scenarios for \texorpdfstring{\cite{salzmann2020trajectron}}{Trajectron++}}\label{sec:appendix-trajectron++}
\begin{figure}[t]
    \centering
    \includegraphics[width=\linewidth]{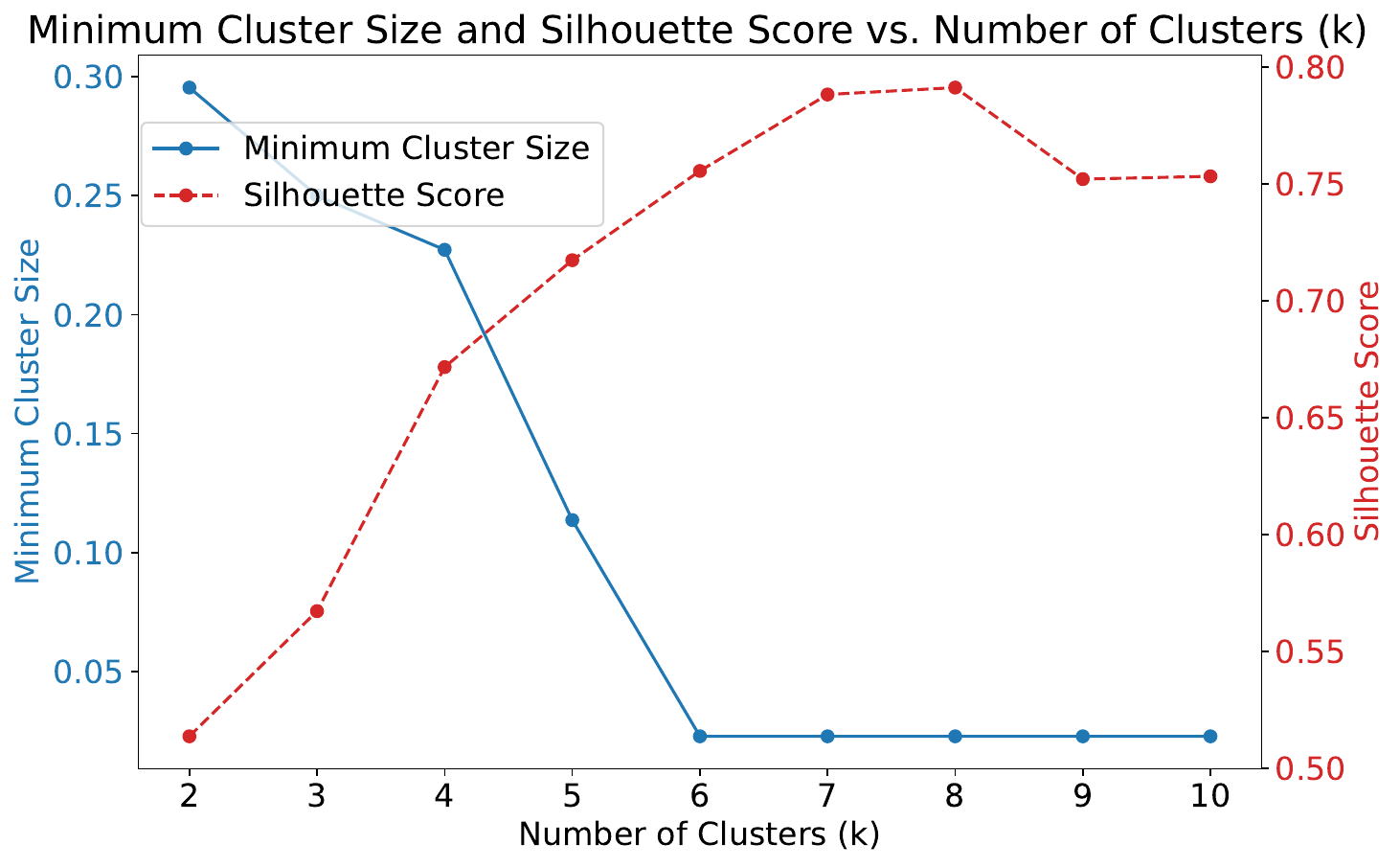}
    \caption{Choice of $k$: a balance of Silhouette Score and minimum cluster size.}
    \label{fig:k-choice-trajectron}
\end{figure}

In this section, we present our findings based on an exemplary behavioral model \cite{salzmann2020trajectron}. The collisions are generated with a prescribed realism parameter $r = 0.03$. The dataset comprises a total of 174 collisions. We have set a threshold for the smallest cluster size at 5\% of the total dataset, aiming to identify a concise set of representative scenarios. For our analysis, we chose $k=5$, illustrated in \Cref{fig:k-choice-trajectron}.

\begin{figure}[b]
    \centering
    \includegraphics[width=0.75\linewidth]{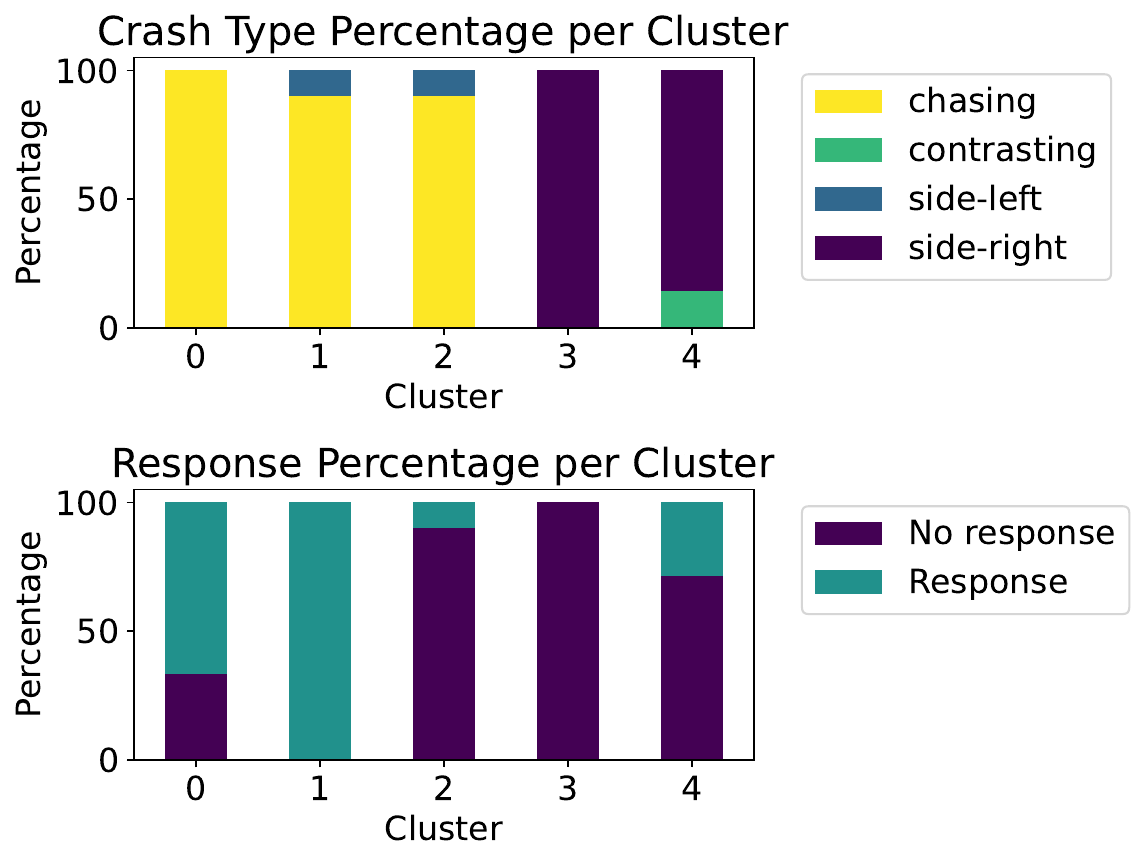}
    \caption{Comparison of angle and response characteristics across resultant scenario clusters.}
    \label{fig:cluster_des_angle_trajectron}
\end{figure}

As illustrated in \Cref{fig:cluster_des_angle_trajectron}, Cluster 1, with relatively higher testing speed, demonstrates an exceptional response rate, with "Chasing" crashes constituting 90\% and "Side-left" crashes 10\%. The complete response rate in this cluster highlights its advanced detection capabilities, showcasing robust performance even under challenging conditions.

In contrast, the other clusters display varying levels of non-responsiveness, each highlighting different challenges in vehicle response systems. Cluster 0: Characterized by the highest testing speed, this cluster exclusively involves "Chasing" crashes (100\%), typically indicating high-speed rear-end collisions. Despite these demanding conditions, the response rate of 67\% reflects that the vehicle systems are generally well-equipped to manage such high-stress scenarios. Cluster 2: With a lower testing speed of 5.80 m/s, this cluster predominantly experiences "Chasing" (90\%) and "Side-left" (10\%) crashes. The notably high no-response rate (90\%) underlines possible inadequacies in the vehicle's sensor accuracy or algorithmic agility, especially in less severe but complex rear-following collisions. Comparing Cluster 0,1,2 together, a deeper investigation into the chasing type collisions is necessary to find out potential system deficiencies.

\begin{figure}[b]
    \centering
    \includegraphics[width=0.75\linewidth]{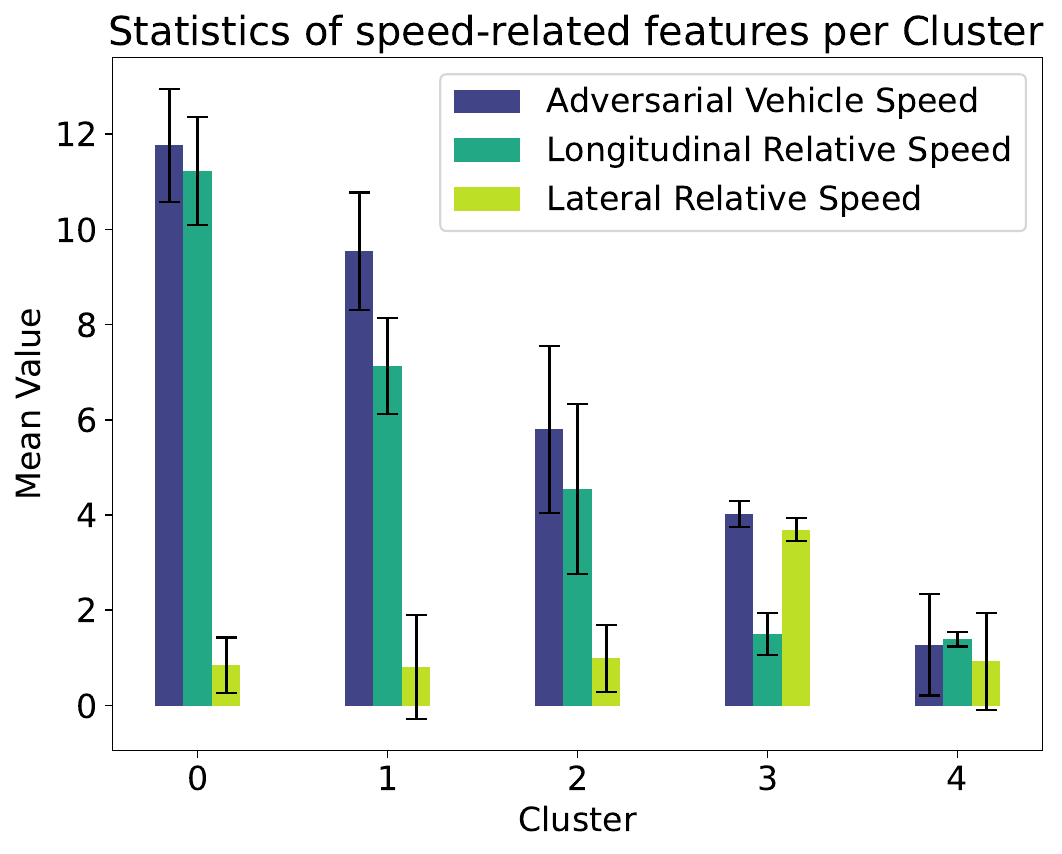}
    \caption{Comparison of clusters via speed-related features.}
    \label{fig:cluster_des_srr_trajectron}
\end{figure}

Cluster 4, characterized by a low testing speed, predominantly involves "Side-right" crashes (86\%) and has a low response rate of 29\%. Cluster 3 exhibits a complete lack of response, with all incidents being "Side-right" crashes at a low speed. These high non-responsiveness rates severely underscore the challenges in detecting lateral threats from the right, pinpointing a critical need for improvements in low-speed collision detection technologies.

Overall, while Cluster 1 sets a benchmark with its full responsiveness, the varied performance across other clusters emphasizes the need for enhanced sensor capabilities and refined algorithms. This disparity particularly necessitates further research into improving detection and response mechanisms in scenarios that currently exhibit high rates of non-responsiveness.

\begin{table*}[h]
    \centering
    \scriptsize
    \begin{tabularx}{\textwidth}{l *{8}{>{\centering\arraybackslash}X}}
        \toprule
        Descriptive Feature & Cluster 0 & Cluster 1 & Cluster 2 & Cluster 3 & Cluster 4 \\
        \midrule
        \textbf{Adversarial Vehicle Speed $v_a$} \\
        Mean & 11.765 & 9.541 & 5.800 & 4.025 & 1.275 \\ 
        Std & 1.183 & 1.234 & 1.761 & 0.270 & 1.063 \\ 
        \midrule
        \textbf{Longitudinal Relative Speed $\Delta v_x$} \\
        Mean & 11.233 & 7.126 & 4.551 & 1.505 & 1.388 \\ 
        Std & 1.132 & 1.007 & 1.779 & 0.442 & 0.147 \\ 
        \midrule
        \textbf{Lateral Relative Speed $\Delta v_y$} \\
        Mean & 0.847 & 0.804 & 0.989 & 3.696 & 0.927 \\ 
        Std & 0.582 & 1.088 & 0.709 & 0.241 & 1.023 \\ 
        \midrule
        \textbf{Crash Type $\beta$} \\
        Chasing & \textbf{48 (100\%)} & \textbf{36 (90\%)} & \textbf{36 (90\%)} & 0 & 0 \\ 
        Contrasting & 0 & 0 & 0 & 0 & 4 (14\%) \\ 
        Side-left & 0 & 4 (10\%) & 4 (10\%) & 0 & 0 \\ 
        Side-right & 0 & 0 & 0 & \textbf{20 (100\%)} & \textbf{24 (86\%)} \\ 
        \midrule
        \midrule
        \textbf{Tested Vehicle Response} \\
        No response & 16 (33\%) & 0 & 36 (90\%) & 20 (100\%) & 20 (71\%) \\
        Response & 32 (67\%) & 40 (100\%) & 4 (10\%) & 0 & 8 (29\%) \\ 
        \bottomrule
    \end{tabularx}
    \caption{Clustering results for the generation based on the behavioral model in \cite{salzmann2020trajectron}. \label{tab:statistics_trajectron}}
\end{table*}


\end{document}